\documentclass[12pt]{amsart}
\usepackage{latexsym, amsmath, amscd, amssymb, amsthm}

\begin{document}

\noindent

 \title{ A note about  the Nowicki conjecture on Weitzenb\"ock derivations }

\author{Leonid Bedratyuk} \address{ Khmelnytsky National University, Instytuts'ka st. 11, Khmelnytsky , 29016, Ukraine}
\email {leonid.uk@gmail.com}

\begin{abstract}
We reduce the Nowicki conjecture on the Weitzenb\"ock derivation of polynomial algebras to well known problem of the classical invariant theory.
\end{abstract}
\maketitle


{\bf 1.} Let $\mathbb{K}$ be a field of characteristic 0. A linear locally
nilpotent derivation $\mathcal{D}$
of the polynomial algebra $\mathbb{K}[Z]=\mathbb{K}[z_1,z_2,\ldots,z_m]$
is called a Weitzenb\"ock derivation.
It is well known that the kernel
$$
\ker \mathcal{D}:=\left\{ f \in \mathbb{K}[Z]| \mathcal{D}(f)=0 \right\}
$$
of the linear locally nilpotent derivation $ \mathcal{D}$ is a finitely
generated algebra,
see \cite{Ses}-\cite{Wei}.

Let $\mathbb{K}[X,Y] = \mathbb{K}[x_1, . . . , x_n,y_1,y_2,\ldots,y_n]$ be the polynomial
$\mathbb{K}$-algebra in $2\,n$ variables. Consider the following Weitzenb\"ock  derivation   $\mathcal{D}_1$  of   $\mathbb{K}[X,Y]:$

$$
\mathcal{D}_1(x_i)=0, \mathcal{D}_1(y_i)=x_i, i=1,2,\ldots, n.
$$

In  \cite{Now} Nowicki conjectured that $\ker \mathcal{D}_1$ is generated by the elements $x_1,$ $x_2,\ldots, x_n$ and the determinants

$$
\left |
\begin{array}{ll}
x_i & x_j \\
y_i &y_j
\end{array}
\right |, 1 \leq i <j \leq n.
$$

The conjecture was confirmed by several authors, see \cite{DrM}, \cite{Kh}, \cite{Kur}.

In this  note we show that the Nowicki conjecture is equivalent to a well known problem of  classical invariant theory, namely, to the problem to  describe  the algebra of joint covariants of $n$ linear binary forms.  Using the same idea we present an explicit set
of generators of the kernel of the derivation $\mathcal{D}_2$ of
$$
\mathbb{K}[X,Y,Z]=\mathbb{K}[x_1,\ldots,x_n,y_1,\ldots,y_n,z_1,\ldots,z_n]
$$
defined by
$$
\mathcal{D}_2(x_i)=0,\mathcal{D}_2(y_i)=x_i,\mathcal{D}_2(z_i)=y_i,\quad
i=1,\ldots,n.
$$

{\bf 2.} It is well known that there is a one-to-one correspondence between  the $\mathbb{G}_a$-actions on an affine algebraic variety $V$   and the locally nilpotent  $\mathbb{K}$-derivations on its algebra of polynomial functions. Let us identify the algebra  $\mathbb{K}[X,Y]$  with the algebra  $\mathcal{O}[\mathbb{K}^{2\,n}]$ of polynomial functions of the algebraic variety   $\mathbb{K}^{2 n}.$    Then, the kernel of the derivation   $\mathcal{D}_1$ coincides with
the invariant ring of the induced  via $\exp(t \, \mathcal{D}_1$  action:
$$
\ker \mathcal{D}_1= \mathcal{O}[\mathbb{K}^{2 n}]^{\mathbb{G}_a}= \mathbb{K}[X,Y]^{\mathbb{G}_a}.
$$

 Now, let $$V_{1}:=\{ \alpha \mathcal{X}+ \beta \mathcal{Y} | \alpha,\beta \in \mathbb{K}\}$$  be the  vector $\mathbb{K}$-space of  linear binary forms endowed with the natural action of the group  $SL_2.$ Consider the induced action of the group  $SL_2$ on the algebra of polynomial functions   $\mathcal{O}[n\,V_{1} \oplus \mathbb{K}^2 ]$   on the vector  space $n\,V_{1} \oplus \mathbb{K}^2,$ where $$n\,V_{1}:=\underbrace{V_{1}\oplus V_{1}\oplus \ldots \oplus V_{1}}_{n \text{ times }} .$$ 
Let  $U_2=\left\{ \left( \begin{array}{cc}  1 & \lambda \\ 0 & 1 \end{array} \right) \Bigl |   \lambda \in \mathbb{K} \right\}$ be the maximal unipotent subgroup of the group  $SL_2.$  The application of the Grosshans principle, see  \cite{Gross}, \cite{Pom}  gives
$$
\mathcal{O}[n\,V_{1} \oplus \mathbb{K}^2 ]^{\,SL_2}\cong \mathcal{O}[n\,V_{1} ]^{\,U_2}.
$$
Thus,
$$
\mathcal{O}[n\,V_{1} \oplus \mathbb{K}^2 ]^{\,\mathfrak{sl_2}}\cong \mathcal{O}[n\,V_{1} ]^{\,\mathfrak{u_2}}.
$$
Since $U_2 \cong (\mathbb{K},+)$   it follows  that $$\ker \mathcal{D}_1 \cong \mathcal{O}[n V_{1} \oplus \mathbb{K}^2 ]^{\,\mathfrak{sl_2}}\cong \mathcal{O}[n\,V_{1} ]^{\,\mathfrak{u_2}}.$$

In the language  of classical invariant theory   the algebra    $C_1:=\mathcal{O}[n V_{1} \oplus \mathbb{K}^2 ]^{\,\mathfrak{sl_2}}$ is called     the  algebra of joint covariants of  $n$  linear binary forms   and the algebra    $S_1:=\mathcal{O}[n V_{1} ]^{\,\mathfrak{u_1}}$ is called     the  algebra of joint semi-invariants of $n$ linear    binary forms. 
Algebras of joint covariants of  binary forms   were an  object of research in the classical invariant theory
of the 19th century.

{\bf 3.}  
 Let us consider the set of  $n$  linear binary forms $f_i=x_i \mathcal{X}+y_i \mathcal{Y},$ $i=1,\ldots,n.$ Then any element of $\mathcal{O}[n V_{1} \oplus \mathbb{K}^2 ]$ can be considered as  a polynomial from  $\mathbb{K}[X,Y,\mathcal{X}, \mathcal{Y}].$
Gordan's famous theorem, see \cite{Gor}, \cite{Gle}, implies:

{\bf Theorem 1}. (A weak form of  Gordan's  theorem)

{\it If $T$ is a subalgebra of $C_1$ with property that $(f_i,z)^r \in T$ whenever $r \in \mathbb{N},$ $z \in T,$  then $T=C_1.$ }

Here $(u,v)^r$ denotes the $r$-transvectants of the binary forms $u,v \in \mathbb{K}[X,Y,\mathcal{X}, \mathcal{Y}]:$
$$
(u,v)^r:=\sum_{i=0}^r (-1)^i { r \choose i } \frac{\partial^r u}{\partial \mathcal{X}^{r-i} \partial \mathcal{Y}^i}   \frac{\partial^r v}{\partial \mathcal{X}^{i} \partial \mathcal{Y}^{r-i}}.
$$

Observe, that $(u,v)^0=u\, v$ and $(u,v)^1$ is exactly  the Jacobian $J(u,v)$ of  the forms $u, v.$ The above theorem yields:

{\bf Theorem 2.} {\it The algebra of joint covariants  $C_1$ of $n$ linear binary forms $f_i,$ $i=1,\ldots,n$  is generated by the forms $f_1, f_2, \ldots, f_n$ and theirs jacobians $J(f_i,f_j),$  $1 \leq i <j \leq n.$
}
\begin{proof}

All forms $f_i,$ $i=1,\ldots,n,$ belong to the algebra of covariants  $C_1.$ By direct calculations we get 
$$
(f_i,f_j)^1=J(f_i,f_j)=\displaystyle \left |
\begin{array}{ll}
\displaystyle \frac{\partial f_i}{\partial \mathcal{X} } & \displaystyle  \frac{\partial f_i}{\partial \mathcal{Y} } \\
& \\
 \displaystyle \frac{\partial f_j}{\partial \mathcal{X} } & \displaystyle  \frac{\partial f_i}{\partial \mathcal{Y} }
\end{array}
\right |=\left |
\begin{array}{ll}
x_i & x_j \\
y_i &y_j
\end{array}
\right |, \text{   and } (f_i,f_j)^r=0, \text{  for }  r>1.
$$
Let us consider the subalgebra  $T$ of $ C_1$ generated by the linears forms $f_1,$ $f_2, \ldots,$ $f_n$ and theirs jacobians $J(f_i,f_j), i<j.$  Since $(f_i, J(f_j,f_k))^r=0$  for all $r \geq 1$ it follows that $T=C_1.$
\end{proof}
Let us show that the result is equivalent to the Nowicki conjecture:

Identify the algebra of semi-invariants $S_1$ with  $\ker \mathcal{D}_1.$ 
The isomorphism  $\tau: C_1 \to S_1 $ take each homogeneous  covariant of degree $m$ (with respect to the variables $\mathcal{X}, \mathcal{Y}$)  to its coefficient of  $\mathcal{X}^m.$ The proof proceeds in the same manner as the proof in the case of unique  binary form, see \cite{Olver}, Proposition 9.45.

Thus, the following statement holds:

{\bf Theorem 3.} {\it The algebra of joint semi-invariants $S_1=\ker \mathcal{D}_1$ of  $n$  linear binary forms $f_i,$ $i=1,\ldots,n$  is generated by the elements  $x_1,x_2, \ldots, x_n$ and theirs determinants 
 $$
\left |
\begin{array}{ll}
x_i & x_j \\
y_i &y_j
\end{array}
\right |, 1 \leq i <j \leq n.
$$
}
\begin{proof} The algebra  $S_1$ is generated by the images of the generating elements  of the algebra $C_1$  under  the homomorphism $\tau.$
We have  $\tau(f_i)=x_i$ and 
$$
\tau(J(f_i,f_j))=\left |
\begin{array}{ll}
x_i & x_j \\
y_i &y_j
\end{array}
\right |.
$$

\end{proof}

 Theorem 3 is exact the Nowicki conjecture.

{\bf 4.} Other ways to prove  Theorem 3 were suggested by  Dersken and Panyushev, see the comments in \cite{DrM}.
Taking into account $ \mathbb{K}^2 \cong_{\mathfrak{sl_2}} V_1$ we get 
$$
\ker \mathcal{D}_1 \cong \mathcal{O}[n V_{1} \oplus \mathbb{K}^2 ]^{\,\mathfrak{sl_2}} \cong \mathcal{O}[(n+1) V_{1}  ]^{\,\mathfrak{sl_2}}.
$$
But then the invariant algebra $ \mathcal{O}[(n+1) V_{1}  ]^{\,\mathfrak{sl_2}}$ is well known because
of the First Fundamental theorem of invariant theory for $SL_2,$ see \cite{Weyl}.

{\bf 5.} A natural generalisation  of the above problem looks as follows.

 Let  $$\mathbb{K}[X,Y,Z] = \mathbb{K}[x_1, . . . , x_n,y_1,y_2,\ldots,y_n,z_1,\ldots,z_n]$$ be the polynomial
$\mathbb{K}$-algebra in $3\,n$ variables. Consider the following derivation  $\mathcal{D}_2$  of the algebra  $\mathbb{K}[X,Y,Z]:$

$$
\mathcal{D}_2(x_i)=0, \mathcal{D}_2(y_i)=x_i, \mathcal{D}_2(z_i)= y_i,i=1,2,\ldots, n.
$$
The following theorem holds:

{\bf Theorem 3.} {\it The kernel of the derivation $\mathcal{D}_2$  is generated by the elements of the following types 
$$
\begin{array}{ll}
1. & x_1,x_2,\ldots, x_n;\\
2. & J_{1,2}, J_{1,3}, \ldots, J_{{n-1},{n}},\\
3. &  H_{1,2}, H_{1,3}, \ldots, H_{{n-1},{n}},\\
4. & \Delta_{1,2,3},\Delta_{1,2,4}, \ldots, \Delta_{{n-2},{n-1},{n}},
\end{array}
$$
where
$$
J_{i,j}:=\left |
\begin{array}{ll}
x_i & x_j \\
y_i &y_j
\end{array}
\right |, 1 \leq i <j \leq n,
$$
$$
H_{i,j}=x_i z_j- y_i y_j+z_i x_j, 1 \leq i \leq j \leq n,
$$
and 
$$
\Delta_{i,j,k}:=\left |
\begin{array}{lll}
x_i & x_j & x_k\\
y_i & y_j & y_k\\
z_i & z_j & z_k\\
\end{array}
\right |, 1 \leq i <j<k \leq n.
$$
}

The proof follows from  the description of generating elements of the algebra of covariants for $n$ quadratic binary forms, see  \cite{GrY}, page 162.

{\bf 5.} 
Any   Weitzenb\"ock derivation of  polynomial algebra is  completely determined by its  Jordan normal form. Denote by   $\mathcal{D}_{k}$   the  Weitzenb\"ok derivation with   Jordan normal form which consists of  $n$  Jordan blocks of size  $k+1.$

{\bf Problem.} Find   a   generating set of $\ker \mathcal{D}_k.$

\end{document}